\documentclass[final]{dmtcs-episciences}

\usepackage[utf8]{inputenc}
\usepackage{subfigure}
\usepackage{graphicx}
\usepackage{amsmath, amssymb, amsthm, tikz, mathrsfs,verbatim,tensor}
\usepackage{setspace}
\usepackage{enumerate}
\usepackage{amsfonts}
\usepackage{latexsym}
\usepackage{graphics}

\usepackage{bbm}

\usepackage{etoolbox}
\patchcmd{\section}{\scshape}{\bfseries}{}{}

\newtheorem{theorem}{Theorem}
\newtheorem{lemma}[theorem]{Lemma}
\newtheorem{corollary}[theorem]{Corollary}

\newtheorem{proposition}[theorem]{Proposition}
\newtheorem{definition}[theorem]{Definition}

\theoremstyle{definition}

\theoremstyle{plain}

\newcommand{\beqlbl}{\begin{equation}}
\newcommand{\eeqlbl}{\end{equation}}

\newcommand{\R}{\ensuremath{\mathbb{R}} }

\newcommand{\N}{\ensuremath{\mathbb{N}} }

\newcommand{\cf}{\mathbf{1}}
\renewcommand{\P}{\mathbb{P}}

\newcommand{\var}{\textrm{Var}}
\newcommand{\E}{\mathbb{E}}

\usepackage{hyperref}
\usepackage{tikz}
\usepackage{xypic} 
\xyoption{curve}

\usetikzlibrary{decorations.markings}
\tikzstyle{vertex}=[circle, draw, inner sep=0pt, minimum size=6pt]
\tikzstyle{Vertex}=[circle, draw, inner sep=0pt, minimum size=14pt]
\tikzstyle{Vertexc}=[circle, draw, inner sep=0pt, minimum size=14pt, fill=blue!30]
\tikzstyle{vertexc}=[circle, draw, inner sep=0pt, minimum size=6pt, fill=red!40]
\tikzstyle{vertexcg}=[circle, draw, inner sep=0pt, minimum size=6pt, fill=green!70!black]

\newcommand{\beq}{\begin{equation*}}
\newcommand{\eeq}{\end{equation*}}
\newcommand{\ba}{\begin{align*}}
\newcommand{\ea}{\end{align*}}

\newcommand{\matbegin}[1]{\left (  \begin{array} {#1} }
\newcommand{\matend}{ \end{array} \right ) }

\newcommand{\fp}{f\!p}

\newcommand{\av}{\textbf{A\!v}}
\newcommand{\iv}{\textbf{Iv}}

\newcommand{\ivn}{\textbf{Iv}_n}

\usepackage[colorinlistoftodos]{todonotes}

\author{Sam Miner\affiliationmark{1}
  \and Douglas Rizzolo\affiliationmark{2}  \and Erik Slivken\affiliationmark{3}\thanks{Partially supported by the ERC Starting Grant 680275 MALIG}
}

\title[Fixed points of pattern-avoiding involutions]{Asymptotic distribution of fixed points of pattern-avoiding involutions
\thanks{The authors would like to the the organizers of Permutation Patterns 2016 at Howard University, during which this project started.}
}
\affiliation{
  Pomona College, Claremont CA, USA\\
  University of Delaware, Newark DE, USA\\
  University of Paris Diderot, Paris, France
}
\received{2017-5-16}
\revised{2017-11-28}
\accepted{2017-12-3}

\begin{document}
	
\publicationdetails{19}{2017}{2}{5}{3658}

\maketitle
\begin{abstract}

For a variety of pattern-avoiding classes, we describe the limiting distribution for the number of fixed points for involutions chosen uniformly at random from that class.  In particular we consider monotone patterns of arbitrary length as well as all patterns of length 3.  For monotone patterns we utilize the connection with standard Young tableaux with at most $k$ rows and involutions avoiding a monotone pattern of length $k$.  For every pattern of length 3 we give the bivariate generating function with respect to fixed points for the involutions that avoid that pattern, and where applicable apply tools from analytic combinatorics to extract information about the limiting distribution from the generating function.  Many well-known distributions appear.

\end{abstract}

%%\author[Christopher~Hoffman]{ \ Christopher~Hoffman$^\star$} 
%\author[Sam~Miner]{ \ Sam~Miner$^{\star}$}
%\author[Douglas~Rizzolo]{ \ Douglas~Rizzolo$^{\dagger}$}
%\author[Erik~Slivken]{ \ Erik~Slivken$^{ \circ}$}
%
%\thanks{\thinspace ${\hspace{-.45ex}}^\star$
%Department of Mathematics,
%Pomona College, Claremont, CA, 91711
%\hfill \\
%\texttt{samuel.miner@gmail.com}
%}
%
%\thanks{\thinspace ${\hspace{-.45ex}}^\dagger$
%Department of Mathematics,
%University of Delaware, Newark, DE, 19716
%\hfill \\
%\texttt{drizzolo@udel.edu}
%}
%
%\thanks{\thinspace ${\hspace{-.45ex}}^\circ$
%LMPA,
%University of Paris Diderot, Paris, France, 75013
%\hfill \\
%\texttt{eslivken@univ-paris-diderot.fr}
%}
%
%
%\date{\today}
%
%\vskip1.3cm
%
%
%
%
%
%
%\maketitle

\section{Introduction}
Identifying the asymptotic distribution of the number of fixed points in a uniformly random permutation is a classic problem in probability whose resolution dates back to Montmort in the early 1700's \cite{de1713essay}, where it is shown that the limiting distribution is Poisson.  Since then, the fixed points of various types of permutations have been intensely studied.  Recently there has been a growing interest in the statistical properties of random pattern-avoiding permutations.  Many of the results concern the overall shape and structure of these permutations \cite{bassino2016brownian, hrs1, madras2015structure, madraspehlivanLD, miner2014shape}, some explore pattern containment \cite{janson2014patterns}, while others consider pattern-avoidance under non uniform distributions such as the Mallows distribution \cite{consecutivepatterns}.  In \cite{hrs3,hrs2} the limiting distribution on the number and location of fixed points is given for a variety of pattern-avoiding classes.  Enumeration for involutions was explored in \cite{simionschmidt} for patterns length 3 and \cite{bona2016pattern} for longer patterns.  In \cite{drs} the number of involutions with a specified number of fixed points was given for each pattern of length 3.  Excellent introductions to the general subject area of pattern avoidance can be found in \cite{bonabook} or \cite{vatter}.

In this paper, we focus on the limiting distribution of the total number of fixed points for random pattern-avoiding involutions.  First we use the machinery of standard Young tableaux to give limiting distributions for involutions with longest increasing or decreasing sequence of length at most $k$.  We then complete the picture for all that avoid a fixed pattern of length $3$.  Taken together, this work suggests that pattern-avoiding permutations are strongly connected to many classical limit theorems in probability.  Indeed, depending on the pattern the (appropriately normalized) limiting distribution of fixed points can be anything from a point mass to a normal distribution to a distribution given in terms of the eigenvalues of the Gaussian Orthogonal Ensemble (conditioned to have trace $0$). 

In some cases, we are able to use straightforward generating function arguments coupled with classical analytic combinatorics to compute the asymptotic distribution of fixed points.  For the sake of the reader, we include an appendix containing the results from analytic combinatorics that we will need.  To our knowledge methods from analytic combinatorics have not been widely used to study asymptotic properties of pattern-avoiding permutations and our results suggest that these methods may be broadly useful for deriving asymptotic permutation statistics.  Some cases, however, cannot be easily done using off-the-shelf results from analytic combinatorics and require different methods.   

We now introduce some necessary notation and summarize our results.  For a permutation $\pi$, let $\fp(\pi)$ be the number of fixed points of $\pi$.  For a permutation $\tau\in S_k$, we say $\pi$  contains the pattern $\tau$ if a there exists a subsequence $i_1<\cdots i_k$ such that $(\pi_{i_1}, \cdots, \pi_{i_k})$ has the same relative order of $\tau$ and $\pi$ avoids $\tau$, or is $\tau$-avoiding, if it does not contain $\tau$. Let $\av_n(\tau)$ be the set of $\tau$-avoiding permutations of $[n]:=\{1,\dots, n\}$.  We are particularly interested in involutions, which are permutations $\pi$ such that $\pi^{-1}=\pi$, and denote the set of $\tau$-avoiding involutions of $[n]$ by $\ivn(\tau)$.   Let $(\xi_{ij})_{1\leq i\leq j\leq k}$ be independent, centered, normal random variables such that $\xi_{ii}$ has variance $1$ and $\xi_{ij}$ has variance $1/2$ if $i<j$.  For $j<i$, define $\xi_{ij}=\xi_{ji}$.  A $k\times k$ random matrix $X$ is said to be drawn from the $k\times k$ Gaussian Orthogonal Ensemble if it is equal in distribution to $(\xi_{ij})_{i,j=1}^k$.  We will need to condition $X$ to have trace $0$.  Since the $(\xi_{ij})_{1\leq i\leq j\leq k}$ are independent normal random variables and the diagonal elements all have the same variance, conditioning the matrix to have trace $0$ is equivalent to projecting $X$ onto the subspace of matrices with trace $0$.  This is part of a larger connection between conditioning normal random variables and projections, see e.g. \cite[Chapter IX]{Janson97}, and is essentially due to the fact that if $I$ is the identity matrix then $X-k^{-1}tr(X)I$ and $k^{-1}tr(X)I$ are independent.  Thus $X$ conditioned to have trace $0$ is equal in distribution to 
\begin{equation}\label{eq goe} 
M = X - \frac{tr(X)}{k}I,
\end{equation}
where $I$ is the identity matrix.  We say that $M$ is a random matrix drawn from the $k\times k$ Gaussian Orthogonal Ensemble conditioned to have trace $0$.  

Our results consider the asymptotic distribution as $n$ increases, and we let $\rightarrow_d$ denote convergence in distribution as $n$ tends to $\infty.$

\begin{theorem}\label{t:321}
Fix $k\in\{2,3,\dots\}$ and let $\Pi_n$ be a uniformly random element of $\ivn((k+1)k\cdots 321)$.  Let $M$ be a random matrix drawn from the $k\times k$ Gaussian Orthogonal Ensemble conditioned to have trace $0$ and let $\Lambda_1\geq \cdots \geq \Lambda_k$ be the ranked eigenvalues of $M$.
\begin{enumerate}
\item If $k$ is even then 
\[\sqrt{\frac{k}{n}} \fp(\Pi_n)  \rightarrow_d \sum_{j =1}^{k} (-1)^{j+1} \Lambda_j.\]
\item If $k$ is odd then 
\[\sqrt{\frac{k}{n}}\left( \fp(\Pi_n)- \frac{n}{k}\right)  \rightarrow_d \sum_{j =1}^{k} (-1)^{j+1} \Lambda_j.\]
\end{enumerate}
\end{theorem}

We remark that when $k+1=3$, we have a simpler description of the limiting distribution.  In this case, we let $B_1$ and $B_2$ be two independent $N(0,1/2)$ random variables (where $N(\mu,\sigma^2)$ denotes a normal random variable with mean $\mu$ and and variance $\sigma^2$), using Equation \eqref{eq goe}, we see that
\[ M =_d \begin{bmatrix}  B_1 &   B_2 \\    B_2 & -B_1 \end{bmatrix}\]
is distributed like a $2\times 2$ GOE matrix conditioned to have trace equal to $0$.  The eigenvalues of this matrix can be computed explicitly, giving 
\[ (\Lambda_1,\Lambda_2) = \left( \sqrt{B_1^2+B_2^2}, -\sqrt{B_1^2+B_2^2}\right),\]
so that if $\Pi_n$ is a uniformly random element of $\iv_n(321)$ then
\[ \sqrt{\frac{2}{n}} \fp(\Pi_n) \rightarrow_d 2 \sqrt{B_1^2+B_2^2} = \sqrt{(2B_1)^2+(2B_2)^2}.\]
Equivalently,
\[ \sqrt{\frac{1}{n}} \fp(\Pi_n) \rightarrow_d  \sqrt{(\sqrt{2}B_1)^2+(\sqrt{2}B_2)^2} .\]
Furthermore, it is well known that $ \sqrt{(\sqrt{2}B_1)^2+(\sqrt{2}B_2)^2}$ follows a Rayleigh($1$) distribution, whose density is given by
\[ f(x) = x e^{-x^2/2} \cf(x\geq 0).\]

\begin{theorem}\label{increasefp}

If $\Pi_n$ is a uniformly random element in $\iv_n(123\cdots k(k+1))$ then

\[
\fp(\Pi_{2n}) \to_d 
 X_{even}
\]
and 
\[
\fp(\Pi_{2n-1}) \to_d 
X_{odd}\,,	
\]

where 
$X_{even}$ has density function given by \[\mathbb{P}(X_{even} = i) = \begin{cases} \frac{\binom{k}{i}}{2^{k-1}} & i \text{ is even}\,,\\
0 & i \text{ is odd}\,,
\end{cases}\]
and $X_{odd}$ has density function given by \[\mathbb{P}(X_{odd} = i) = \begin{cases} \frac{\binom{k}{i}}{2^{k-1}} & i \text{ is odd}\,,\\
0 & i \text{ is even}\,.
\end{cases}\]
  	
\end{theorem}

\begin{theorem}\label{t:231}
Fix $\tau \in \{ 231, 312\}$.  If $\Pi_n$ is a uniformly random element of $\iv_n(\tau)$ then $\E(\fp(\Pi_n))=\frac n3 + O(1)$, $\var(\fp(\Pi_n)) = \frac{8}{27}n + O(1)$ and 
\[ \frac{\fp(\Pi_n) - \frac{1}{3}n}{ \sqrt{8n/27}} \rightarrow_d Z,\]
where $Z$ is a standard normal random variable.
 \end{theorem}

To the best of our knowledge, the above theorems are new.  The following theorem was first established in \cite{hrs3} using the theory of local limits for Galton-Watson trees as a corollary of a stronger result that also described the location of the fixed points.  We include a simple proof establishing the asymptotic distribution of the total number of fixed points using analytic combinatorics.

\begin{theorem}[Corollary 6 \cite{hrs3}]\label{avn321}
Fix $\tau \in \{321, 213, 132\}$.  If $\Pi_n$ is a uniformly random element of $\av_n(\tau)$, then
\[ \fp(\Pi_n) \rightarrow_d N \]
where $N$ has a negative-binomial distribution with parameters $r=2$ and $p=1/3$ ( $\P(N = k) = \frac {4(k+1)} {9} (1/3)^k $).
\end{theorem}

\proof[ of Theorem \ref{avn321}]

We begin with a bivariate generating function from Elizalde \cite{elizaldethesis} for $\tau = 321, 213, 132$, 

\begin{equation}\label{eli}
F_{\av(\tau)}(x,t)=\sum_{\sigma\in \av(\tau)} x^{\fp(\sigma)}t^{|\sigma|} =\frac{2}{1+2t(1-x)+\sqrt{1-4t}}.
\end{equation}

Note that for fixed $x$, with say $|x|< 2$, we have that 
\[ \lim_{t\to 1/4}  \frac{1}{\sqrt{1-4t}} \left(  \frac{ 2}{1+ 2t( 1-x)+ \sqrt{1-4t} } - \frac{4}{3-x} + \frac{8}{(3-x)^2} \sqrt{1-4t}\right) =0 .\]
It follows from Corollary \ref{VI.1}, which is \cite[Corollary VI.1]{FlSe09}, that 
\[ [t^n]F_{\av(\tau)}(x,t) \sim -[t^n]\frac{8}{(3-x)^2} \sqrt{1-4t} \sim \frac{4^{n+1}}{(3-x)^2\sqrt{\pi n^3} }.\]
If $\sigma_n\in \av_n(\tau)$ is uniformly random, we have that
\[ f_n(x) = \sum_{k=0}^\infty \P(\fp(\sigma_n)=k) x^k = \frac{[t^n]F_{\av(\tau)}(x,t)}{[t^n]F_{\av(\tau)}(1,t)} \rightarrow   \frac{4}{(3-x)^2} = \sum_{k=0}^\infty \frac{4}{9}(k+1) \left(\frac{1}{3}\right)^k x^k .\]
The result now follows from the standard continuity result for probability generating functions, see e.g. Theorem \ref{IX.1}, which is \cite[Theorem IX.1]{FlSe09}. 
\qed

The next two theorems, which have already been established in \cite{hrs2}, have significantly more complicated proofs than the previous theorems we have stated. In fact, \cite[Theorem 1.1]{hrs2} gives full information not only about the total number of fixed points, but also their locations.  We do not give proofs of them, but include them here in order to give a complete catalogue of results for avoidance of patterns of length $3$.

\begin{theorem}[Theorem 1.1 \cite{hrs2}]
Let $(\mathbbm{e}_t,0\leq t\leq 1)$ be a standard Brownian excursion and let $\Pi_n$ be a uniformly random element of $\av_n(231)$ or $\av_n(312)$.  Then
\[ \frac{1}{n^{1/4}} \fp(\Pi_n) \rightarrow_d \frac{1}{2^{7/4}\sqrt{\pi}}\int_0^1 \frac{1}{\mathbbm{e}_t^{3/2}} dt.\]
\end{theorem}

\begin{theorem}[Theorem 1.1 \cite{hrs2}]
Let $\Pi_n$ be a uniformly random element of $\av_n(123)$ and let $A$ and $B$ be independent Bernoulli($1/4$) random variables.  Then
\[ \fp(\Pi_n) \rightarrow_d A+B.\]
\end{theorem}

\section{Involutions avoiding large monotone patterns}

Before we proceed we give a few basic definitions concerning Young diagrams and Young tableaux.  A Young diagram is a collection of rows of boxes left justified so that the number of boxes in each row is weakly decreasing.  A Young tableau is a Young diagram with each box filled with a number from $\N$ such that the numbers are weakly increasing along rows and columns.  A standard Young tableau of size $n$ is a Young tableau of size $n$ where each element in $[n]$ appears exactly once.  Given a Young tableau, $T$, we call the corresponding Young diagram the shape of $T$ and denote it $\lambda(T)$.  We denote the conjugate of a Young diagram $\lambda$ by $\lambda'$, where the row counts of $\lambda'$ are given by the column counts of $\lambda$.  The number of boxes in the $i$th row is denoted by $\lambda_i$, and hence the number of boxes in the $i$th column is denoted by $\lambda'_i$.  The conjugate of $T$, denoted $T'$, is also a standard Young tableaux with shape $\lambda(T') = \lambda'(T).$

Denote the set of Young diagrams of size $n$ with at most $k$ rows by $\mathcal{D}_{n,k}$ and the set of Young tableaux of size $n$ with at most $k$ rows by $\mathcal{T}_{n,k}.$  For a fixed $\nu \in \mathcal{D}_{n,k}$ we let $\mathcal{T}_{n,k}(\nu)$ denote the subset of $\mathcal{T}_{n,k}$ with shape $\nu$.  For $T \in \mathcal{T}_{n,k}$, and $m\leq n$, we let $T^{m}$ denote the standard Young tableau contained in $T$ with entries $1,\cdots,m.$  For $\gamma\in \mathcal{D}_{m,k}$ we let $\mathcal{T}_{n,k}(\gamma)$ denote the subset of $\mathcal{T}_{n,k}$ such that $\lambda(T^m) = \gamma.$  Similarly we let $\mathcal{D}_{n,k}(\gamma)$ the subset of $\mathcal{D}_{n,k}$ that contain $\gamma$.  Finally we let $\mathcal{T}_{n,k}(\nu,\gamma) = \mathcal{T}_{n,k}(\nu) \cap \mathcal{T}_{n,k}(\gamma).$  

For $\nu \in \mathcal{D}_{n,k}(\gamma)$ we let $\nu/\gamma$ denote the skew Young diagram obtained by removing the boxes of $\gamma$ from $\nu$.  If $T\in \mathcal{T}_{n,k}(\gamma,\nu)$ we construct the skew standard Young tableaux of shape $\nu/\gamma$ by removing from $T$ the boxes associated with $T^m$.  We denote this skew standard Young tableau by $T^{m\nearrow n}$, the set of skew standard Young tableaux of shape $\nu/\gamma$ by $\mathcal{S}_{n,k}(\nu/\gamma)$, and the union over $\nu$ of these sets by $\mathcal{S}_{n,k}(\gamma).$  The set $ \mathcal{T}_{n,k}(\gamma)$ is in bijection with the direct sum $\mathcal{T}_{m,k}(\gamma)\oplus \mathcal{S}_{n,k}(\gamma).$

The Robinson-Schensted-Knuth (RSK) algorithm gives a bijection between pairs of standard Young tableaux $(P,Q)$ and permutations of length $n$.  If $\pi$ is an involution then the corresponding pair satisfies $P=Q$.  For an involution $\pi$ we let $T=T(\pi)$ denote the unique standard Young tableaux obtained by RSK. The number of rows in $T$ is giving by the longest decreasing sequence in $\pi$.  Similarly the number columns of $T$ is given by  given by the longest increasing sequence \cite{schensted}.  Let $\tau_k$ denote the monotone decreasing pattern $(k+1)k,\cdots,321$ and $\rho_k$ the reverse of $\tau_k$.  RSK gives a bijection between $\mathcal{T}_{n,k}$ and $\iv_n(\tau_k)$.  By conjugation $\mathcal{T}_{n,k}$ is also in bijection with $\iv_n(\rho_k)$.  For a modern reference of the RSK algorithm see \cite{stanley2}

The number of fixed points of an involution is equal to the number of odd columns of the corresponding tableau \cite{schutzenberger}.  The following proposition follows from results found in \cite{Mat08} and in \cite{Sn07} that consider Young tableaux with bounded number of rows or columns.  

\begin{proposition}\label{gaps}
Let $\tau_k$ denote the monotone decreasing pattern $(k+1)k,\cdots,321$ and let $(\Lambda_i)_{1\leq i\leq k}$ be the ranked eigenvalues of a traceless GOE matrix.  For $\Pi_n$ chosen uniformly from $\iv_n(\tau_k)$, 

\[\left(\sqrt{\frac{k}{n}}\left( \lambda_i(\Pi_n)-\frac{n}{k}\right)\right)_{1\leq i \leq k} \longrightarrow_d \  (\Lambda_i)_{1\leq i\leq k}.
\]

Moreover, for any fixed $d > 0,$

\[\P\left( \min_{2\leq i \leq k}\{ \lambda_{i-1}(\Pi_n) -\lambda_{i}(\Pi_n) \} < d \right) \to 0.\]
   	
\end{proposition}

\proof Letting $\Pi_n$ be a uniformly random element of $\iv(\tau_k)$ it follows from the proof of \cite[Theorem 3]{Sn07} by taking the square root of Equation (10) there or by \cite[Theorem 5.1]{Mat08} by the implied local limit theorem that
\[ \left( \sqrt{\frac{k}{n}}\left( \lambda_i(\Pi_n) - \frac{n}{k}\right)\right)_{1\leq i\leq k}\]
converge in joint distribution to the vector $(\Lambda_i)_{1\leq i\leq k}$ of ranked eigenvalues of a traceless GOE matrix.  If we let 
\[ \mathfrak{H}_k = \left\{ (x_1,\dots, x_k) \in \R^k \middle | x_1\geq \cdots\geq x_k, \ x_1+\cdots + x_k=0\right\},\]
then the probability density function of $(\Lambda_i)_{1\leq i\leq k}$ with respect to the $(k-1)$-dimensional Hausdorff measure on $\mathfrak{H}_k$ is
\[g(x_1,\cdots,x_k)= \frac{1}{Z_k} e^{-\frac{1}{2} \sum_{j=1}^k x_j^2} \prod_{1\leq i<j \leq k} (x_i-x_j),\]
where $Z_k$ is a normalizing constant.  For each $2\leq i \leq k$, the function $h_i =g/(x_{i-1}-x_i)$ is integrable on $\mathfrak{h}_k$ and satisfies $g < \epsilon h_i$ on the subset $\{x_{i-1} -x_i < \epsilon\}$. By the integrability of $h$ over $\mathfrak{h}_k$, there exists $M>0$ such that, along with the union bound
\[ 
\P\left( \min_{2\leq i \leq k}\{\Lambda_{i-1} - \Lambda_i\} < \epsilon \right) \leq \sum_i\int_{\mathfrak{h}_k } \epsilon h_i \leq \epsilon kM.
\]

For any $\epsilon >0$ and large enough $n$, $d\sqrt{\frac{k}{n}} < \epsilon$ finishing the proof.
\qed
%\[\begin{split}\P(\lambda_{i-1}(\Pi_n)-\lambda_i(\Pi_n) < d) &= \P( \sqrt{\frac{k}{n}}\lambda_{i-1}(\Pi_n)-\lambda_i(\Pi_n) < \sqrt{\frac{k}{n}} d ) \\
%&\to \P( \Lambda_{i-1}-\Lambda_i \leq \epsilon ). 	
%\end{split}
%\]

\proof[of Theorem \ref{t:321}]

If $\pi\in  \iv(\tau_k)$ then $\lambda_j(\pi) = 0$ for $j>k$ and the number of odd columns (hence the number of fixed points \cite{schutzenberger}) is given by 
\[ \fp(\pi) = \sum_{j =1}^{k} (-1)^{j+1} \lambda_j(\pi) = \sum_{j =1 \atop j \textrm{ odd}}^{k}  \lambda_j(\pi) - \sum_{j =1 \atop j \textrm{ even}}^{k}  \lambda_j(\pi).\]

Consequently by Proposition \ref{gaps}, if $k+1$ is odd, so $k$ is even, we have
\[\begin{split} \sqrt{\frac{k}{n}} \fp(\Pi_n) =  \sqrt{\frac{k}{n}} \left( \fp(\Pi_n)- \sum_{j=1}^k (-1)^{j} \frac{n}{k}\right) & =   \sum_{j =1}^{k} (-1)^{j+1} \sqrt{\frac{k}{n}} \left(\lambda_j(\Pi_n) - \frac{n}{k}\right)\\
& \rightarrow_d  \sum_{j =1}^{k} (-1)^{j+1} \Lambda_j, \end{split} \]
while if $k+1$ is even, so $k$ is odd, then 
\[\begin{split} \sqrt{\frac{k}{n}}\left( \fp(\Pi_n)- \frac{n}{k}\right) =  \sqrt{\frac{k}{n}} \left( \fp(\Pi_n)- \sum_{j=1}^k (-1)^{j} \frac{n}{k}\right) & =   \sum_{j =1}^{k} (-1)^{j+1} \sqrt{\frac{k}{n}} \left(\lambda_j(\Pi_n) - \frac{n}{k}\right)\\
& \rightarrow_d \sum_{j =1}^{k} (-1)^{j+1} \Lambda_j. \end{split} 
\]
\qed

Before we begin the proof of Theorem \ref{increasefp}, we define a Markov chain whose stationary distribution will precisely describe the limiting distribution of the number of fixed points. 

Consider the Markov chain $C$ with state space $S = \{0,1,\ldots,k\}$, and transition matrix $P$ with probabilities 
\[P_{i,j} = \begin{cases} \frac{i}{k} & j = i-1\\
1-\frac{i}{k} & j = i+1\\
0 & \text{ otherwise.}
\end{cases}
\]
This is a discrete version of the Ehrenfest urn model~\cite{Ehren}, which can be interpreted as having $k$ balls divided between two urns, and at each step choosing a ball uniformly, and moving it to the other urn. The number of rows of odd length in a Young diagram with at most $k$ rows will have the same state space as this Markov chain.  Adding a box to a row will change its parity and therefore is equivalent to moving a ball from one urn to the other.  The difficulty with Young diagrams is that they do not grow by choosing a row uniformly at random and adding a box to it, since the rows must be decreasing in order.  Dealing with this complication is the main difficulty in the proof of Theorem \ref{increasefp}.

Known results about the Ehrenfest model give us the following lemma.

\begin{lemma}\label{ehrenfest} The chain $C$ is periodic of period 2.   If the initial position of the chain is deterministic then as $d \to \infty$, $C_d$ approaches alternation between vectors ${\bf p}$ and ${\bf q} \in S$, where 
\[{\bf p}_i = \begin{cases} \frac{\binom{k}{i}}{2^{k-1}} & i \text{ is even}\\
0 & \text{ otherwise,}
\end{cases}
\]
and
\[{\bf q}_i = \begin{cases} \frac{\binom{k}{i}}{2^{k-1}} & i \text{ is odd}\\
0 & \text{ otherwise.}
\end{cases}
\]

\end{lemma}

\proof[of Lemma~\ref{ehrenfest}]
First, $C$ is irreducible since every state can reach every other state, and therefore positive recurrent (since the state space is finite). The periodicity of $C$ is clearly 2, since at each transition the parity of our state changes. Claim: the invariant probability density function of $C$ is given by \[f(i) = \binom{k}{i}\left(\frac{1}{2}\right)^k\,, \qquad i \in S\,.\]
This is true since applying one transition to $f$ yields 
\[
\aligned
fP(i) &= f(i-1)p_{i-1,i} + f(i+1)p_{i+1,i}\\
&= \binom{k}{i-1}\left(\frac{1}{2}\right)^k\left(1-\frac{i-1}{k}\right) + \binom{k}{i+1}\left(\frac{1}{2}\right)^k\frac{i+1}{k}\\
&= \left(\frac{1}{2}\right)^k \left[\binom{k-1}{i-1}+\binom{k-1}{i}\right]\\
&= \left(\frac{1}{2}\right)^k \binom{k}{i} = f(i)\,.
\endaligned
\]
Applying the periodicity of $C$ to this invariant density function completes the proof.
\qed

\proof[ of Theorem \ref{increasefp}]

For $\pi\in \iv_n(\rho_k)$ the number of fixed points of $\pi$ is at most $k$ and the parity of the number of fixed points must match the parity of $n$.  The number of fixed points will be given by the number of rows of odd length in the conjugate of $T(\pi).$   

For $\lambda\in \mathcal{D}_{n,k}$ let $f_\lambda$ denote the size of $\mathcal{T}_{n,k}(\lambda)$.  For $m\leq n$ and $\gamma\in \mathcal{D}_{m,k}$ let $s_\gamma$ denote the size of $\mathcal{S}_{n,k}(\gamma)$.  The decomposition of $T\in \mathcal{T}_{n,k}$ into a direct sum implies
\[
	|\mathcal{T}_{n,k}(\gamma)| = f_\gamma s_\gamma.
\]
For distinct $\gamma, \mu \in \mathcal{D}_{m,k}$ the sets $\mathcal{T}_{n,k}(\gamma)$ and $\mathcal{T}_{n,k}(\mu)$ are disjoint, so
\[
|\mathcal{T}_{n,k}| = \sum_{\gamma\in \mathcal{D}_{m,k}} f_\gamma s_\gamma. 
\]

We wish to consider a random standard Young tableau	 chosen uniformly from $\mathcal{T}_{n,k}.$  We construct the following probability measure $w_1$ on $\mathcal{D}_{m,k}$ where 
\[
w_1(\gamma) := \frac{f_\gamma s_\gamma}{\sum_{\mu\in \mathcal{D}_{m,k}}f_\mu s_\mu}.
\]

Conditioned on $\gamma$, choosing $T^m\in \mathcal{T}_{m,k}(\gamma)$ uniformly with probability $f_\gamma^{-1}$ and $T^{m\nearrow n}\in \mathcal{S}_{n,k}(\gamma) $ with probability $s_\gamma^{-1}$ gives a uniform random standard Young tableau $T=T^m \oplus T^{m\nearrow n}$ in $\mathcal{T}_{n,k}(\gamma)$.  By choosing $\gamma$ with probability $w_1(\gamma)$ this process gives a uniformly random element of $\mathcal{T}_{n,k}$.   

Fix $d>0$ and let $m = n-d.$  Let $\delta(\gamma) = \min_{2\leq i \leq k} (\gamma_{i-1} - \gamma_i)$ denote the minimum difference in the lengths of consecutive rows.  If $\delta(\gamma) > d$ then every skew standard Young tableaux in $\mathcal{S}_{n,k}(\gamma)$ will consist of up to $k$ non-overlapping horizontal strips.  For every such $\gamma$, we have $s_\gamma = k^{d}$.  For $\gamma$ such that $\delta(\gamma) < d$ we have $s_\gamma < k^d.$

We define another measure on $\mathcal{D}_{n,k}$,
\[
w_2(\gamma) :=  \frac{k^df_\gamma}{ \sum_\gamma k^df_\gamma} = \frac{f_\gamma}{\sum_{\gamma}{f_\gamma}}.
\]
The measure $w_2$ is equivalent to choosing $T^m$ uniformly from $\mathcal{T}_{m,k}$ and letting $\gamma = \lambda(T^m).$  Then 
\[w_1(\delta(\gamma) < d) \leq w_2( \delta(\gamma) < d).\]
By Proposition \ref{gaps}, as $m$ increases, $\P_{w_2}(\delta(\gamma))<d)\to 0$, and therefore $\P_{w_1}(\delta(\gamma )< d )\to 0$.

Let $Y$ denote the number of odd rows of $\gamma$ chosen with probability $w_1(\gamma)$.  Both $m$ and $Y$ have the same parity.  With probability tending to 1 as $m$ increases, $\delta(\gamma) >d$.  Conditioned on $\delta(\gamma)> d$, we may construct $T^{m\nearrow n}\in \mathcal{S}_{n,k}(\gamma)$ by choosing uniformly from each of the $k$ rows and a adding box to that row, repeating this process $d$ times.    The conditional distribution of the number of odd rows of $T = T^m \oplus T^{m \nearrow n}$ on is given exactly by the Markov chain in Lemma \ref{ehrenfest}.  If $m$ and $d$ have the same parity (hence $n=m+d$ is even), then $C_{d}(Y)\to_d X_{even}$.  Otherwise $C_{d}(Y) \to_d X_{odd}$.   
\qed

% = \P\left(\min_{2\leq i\leq k}\left( \sqrt{\frac{k}{2m}}( \lambda_{i-1}(\Pi_{2m})  - \lambda_{i}(\Pi_{2m}) \right) < O\left(\epsilon\sqrt{k} m^{-0.4}\right) \right) \to 0.$$

\section{Bivariate Generating Functions}\label{bgfs}

Bivariate generating functions are valuable tools in understanding limiting statistics of combinatorial classes.  For a fixed $\tau\in S_3$ we define the following bivariate generating function with respect to fixed points $\fp$:
\begin{equation}
	F_{\iv(\tau)}(x,t) = \sum_{\sigma \in \iv(\tau)} x^{\fp(\sigma)}t^{|\sigma|}.
\end{equation}
  There are three distinct generating functions.  One for $\tau \in \{231, 312\}$, one for $\tau \in \{321, 132, 213\},$ and another for $\tau=123$. For our purposes we will only use $F_{\iv(231)}(x,t),$ though we include the others for completeness.

\begin{proposition}\label{bgfs231}

For $\iv(231) = \iv(312)$ we have 
\[F_{\iv_n(231)}(x,t)= \sum_{\sigma\in \iv(231)} x^{\fp(\sigma)}t^{|\sigma|} = \frac{1-t^2}{1-2t^2-xt}.
\]
%\frac{1-t^2}{1-2t^2-t^2\left(\frac{1+xt}{1-xt}\right)}.$$
\end{proposition}
\proof

If a permutation $\pi$ avoids the pattern $\tau$ then $\pi^{-1}$ avoids the pattern $\tau^{-1}$.  Therefore involutions that avoid the pattern 231 must also avoid 312.  This allows for a bijection between $\iv(231)$ and the set of compositions of integers into positive parts \cite[Proposition 6]{simionschmidt}.  For example 

\beqlbl\label{ex1}
3 \ 2 \ 1  \ 5 \ 4 \ 8 \ 7 \ 6 \longleftrightarrow (3,2,3).
\eeqlbl

If $(a_1, \cdots, a_k)$ is a composition, $b_1=0$, and $b_i = \sum_{j=1}^{i-1}i a_j$, then the corresponding permutation, $\pi$, under the bijection is defined point-wise for $b_i < t \leq b_{i+1}$ by 
\[
\pi(t) := b_i + a_i - (t-b_i) +1.
\]
A fixed point occurs if and only if $t = (a_i + 1)/2+bi$ which can only occur if $a_i$ is odd and can only occur once for each odd $a_i$.  Hence, there is one fixed point for every part of odd size.  In Example \ref{ex1} both $2$ and $7$ are fixed points and occur in the middle of an odd decreasing sequence.  

Let $A(n,k)$ denote the set of compositions of $n$ into exactly $k$ positive parts and $A_{odd}(n,k)$ be the set of compositions of $n$ into exactly $k$ positive odd parts.  For each $\rho \in A(n,k)$ with exactly $k-j$ odd parts there is a unique $\nu\in A_{odd}(n-j,k)$ obtained by reducing each even part by one.  For each $\nu\in A_{odd}(n,k)$ and each $0\leq j\leq k$ there are precisely ${k \choose j}$ permutations in $A(n+j,k)$ with exactly $k-j$ odd parts that reduce to $\nu.$

 We then have the following:

\begin{align*}
\sum_{\sigma\in \iv(231)}x^{\fp(\sigma)}t^{|\sigma|} &= \sum_{n,k}\sum_{\rho\in A(n,k) }t^n x^{\# \{\text{ odd parts of $\rho$ }\}}\\
&= \sum_{n,k,j}\sum_{\nu \in A_{odd}(n,k)}{k\choose j} t^{n+j}x^{k-j}\\
&= \sum_{j,k}\sum_{s_1, \cdots ,s_k }{k\choose j}(t^2)^{s_1+\cdots +s_k}t^{k+j}x^{k-j}\\
%&= \sum_{j,k}{k\choose j}\left(\frac tx\right)^j \left(\frac{xt}{1-t^2}\right)^k\\
%&= \sum_k\left( \frac{xt+t^2 }{1-t^2} \right)^k\\
&= \frac{1-t^2}{1-2t^2-xt}
\end{align*}
\qed

\proof[of Theorem~\ref{t:231}]

We now set up the notation to apply Theorem \ref{IX.9}, which \cite[Theorem IX.9]{FlSe09}.  Let $B(x,t) = 1-t^2$ and $C(x,t) = 1-2t^2-xt$, so that $F_{\iv(231)}(x,t) = B(x,t)/C(x,t)$, with $B$ and $C$ being analytic.  Since $C(1,t)=1-2t^2-t=-2(t+1)(t-1/2)$ we see that that $F_{\iv(231)}(1,t) $ is meromorphic on the ball of radius $r=3/4$ with only a simple pole at $\rho=1/2<3/4$.  Moreover, $B(1,\rho)\neq 0$.  Additionally, 
\[ C_t(1,\rho)C_x(1,\rho) = (-4\rho-1)(-\rho) = 3/2\neq 0,\]
and if we let $\rho(x) = (-x +\sqrt{x^2+8})/4$ then $\rho(1)=\rho$ and $C(x,\rho(x))=0$ and $\rho$ is analytic at $1$.  Let $f(x) = \rho(1)/\rho(x) = 1/(2\rho(x))$.  Note that
\[\mathfrak{v}\left(\frac{\rho(1)}{\rho(x)}\right)=\mathfrak{v}\left(f(x)\right) :=  \frac{f''(1)}{f(1)}+\frac{f'(1)}{f(1)} - \left(\frac{f'(1)}{f(1)}\right)^2 = \frac{2}{27} + \frac{1}{3} - \frac{1}{9} = 8/27,\]
where we have used that $f(1)=1$, $f'(1)=1/3$, and $f''(1)=2/27$.  Similarly, we see that
\[ \mathfrak{m}\left(\frac{\rho(1)}{\rho(x)}\right)=\mathfrak{m}(f(x)) := \frac{f'(1)}{f(1)}=\frac{1}{3}.\]
If we let $X_n$ be the number of fixed points of a uniformly random $231$-avoiding involution of $[n]$, then $\P(X_n=k) = [x^kt^n]F_{\iv(231)}(x,t)/[t^n]F_{\iv(231)}(1,t)$, and it follows from Theorem IX.9 of \cite{FlSe09} that $\E(X_n) = (n/3)+O(1)$, $Var(X_n) = (8n/27)+O(1)$, and that
\[ \P\left(  \frac{X_n - \E(X_n)}{ \sqrt{Var(X_n)}} \leq x\right) = \P(Z\leq x) + O(1/\sqrt{n}),\]
where $Z$ is a standard normal random variable.  It follows immediately that
\[ \frac{X_n - \frac{1}{3}n}{ \sqrt{8n/27}} \rightarrow_d Z. \qedhere\]
\qed\newpage

\begin{proposition}
	
Fix $\tau \in \{ 321, 132, 213\}.$  For $\iv(\tau)$ we have
\[F_{\iv(\tau)}(x,t)   = \frac{ 2}{1-2xt+\sqrt{1-4t^2}}.\]

\end{proposition}

\proof

We first consider $F_{\iv(321)}(x,t)$.  Permutations in $\ivn(321)$ are in bijection with symmetric Dyck paths of length $2n$ and under this bijection the centered tunnels correspond to fixed points \cite{elizalde2004bijections}.  Let $D(n,k)$ denote the set of symmetric Dyck paths of length $2n$ with exactly $k$ centered tunnels.  Any path in $D(n,k)$ is uniquely determined by its first $n$ steps.  Let $u$ and $d$ denote the number of up and down steps respectively in the first half of Dyck path.  For paths in $D(n,k)$ we have $u+d=n$ and $u-d=k$.  By standard ballot counting arguments we have 
\[|D(2d+k,k)| = {2d+k\choose d} - {2d+k \choose d-1}\]
so that 
\begin{equation}\label{fpequation321}
\sum_{\pi\in\iv(321)}t^nx^{\fp(x)} = \sum_{d,k\geq 0}\left({2d+k\choose d} - {2d+k \choose d-1}\right) t^{2d+k}x^k.
\end{equation}

Using the identity $\sum_{m\geq 0} {2n+m\choose n}s^n = \frac{1}{\sqrt{1-4s}}\left(\frac{1-\sqrt{1-4s}}{2s} \right)^m$ we have
\begin{equation}\label{part1}
\sum_{d\geq 0}{2d+k\choose d}t^{2d+k}x^k = \frac{(xt)^k}{\sqrt{1-4t^2}}\left( \frac{(1-\sqrt{1-4t^2})}{2t^2} \right )^k
\end{equation}
and 
\begin{equation}\label{part2}
	\sum_{d\geq 0}{2(d-1)+k+2\choose d-1}t^{2(d-1)+k+2}x^k = \frac{(xt)^kt^2}{\sqrt{1-4t^2}}\left( \frac{(1-\sqrt{1-4t^2})}{2t^2} \right )^{k+2}
\end{equation}

Combining \eqref{part1} and \eqref{part2} with the appropriate factors and summing over $k\geq 0$ gives

\[\sum_{d,k\geq 0}\left({2d+k\choose d} - {2d+k \choose d-1}\right) t^{2d+k}x^k=\frac{1-t^2\left(\frac{1-\sqrt{1-4t^2}}{2t^2}\right)^2}{\sqrt{1-4t^2}}\left(\sum_{k\geq 0}(xt)^k\left(\frac{1-\sqrt{1-4t^2} }{2t^2}\right) ^k \right)
\]
which allows us to simplify \eqref{fpequation321} to
\[
	\sum_{\pi\in\iv(321)}t^nx^{\fp(x)}=\frac{2}{1-2xt+ \sqrt{1-4t^2}}
\] with some straightforward manipulations.

Rotation by 180 degrees sends $\ivn(132)$ to $\ivn(213)$ and preserves the number of fixed points.	 From \cite[Theorem 8]{elizalde2004bijections}, or \cite[Theorem 2.3]{drs} there is bijection from $\av_n(321)$ to $\av_n(132)$ which preserves the number of fixed points and also commutes with taking inverses, so induces a bijection between $\ivn(321)$ and $\ivn(132).$  Therefore \[F_{\iv(321)}(x,t) = F_{\iv(132)}(x,t) = F_{\iv(213)}(x,t).\]
\qed

\begin{corollary}
Fix $\tau\in \{321,132,213\}$.  If $\Pi_n$ is a uniformly random element of $\iv_n(\tau)$, then
\[
\sqrt{\frac{1}{n}}\fp(\Pi_n)\to_d X,
\]
where $X$ follows a Rayleigh(1) distribution.
	
\end{corollary}

\proof This follows from Theorem \ref{t:321} and  \cite[Theorem 8]{elizalde2004bijections} or \cite[Theorem 2.3]{drs}.
\qed

For completeness we also give the bivariate generating function for $\iv(123)$ with respect to fixed points.  

\begin{proposition}

For $\iv(123)$ we have

\[F_{\iv(123)}(x,t)= \sum_{\rho \in \iv{123}} x^{\fp(\rho)}t^{|\rho |} = 1+(tx+t^2(1+x^2))\left(\frac{1-\sqrt{1-4t^2}}{2t^2\sqrt{1-4t^2}}\right)\,.\]

\end{proposition}

\proof
By the RSK correspondence, permutations in $\ivn(123)$ are in one-to-one correspondence with standard Young tableaux of size $n$ with at most two columns. The number of fixed points in such a permutation is equal to the number of columns of odd length in the standard Young tableaux. If $n$ is odd, there is exactly one odd column in the Young tableaux. If $n$ is even (so $n=2k$), there are either zero or two columns. Every tableaux of size $2k$ can be created by placing the element $2k$ in either the first or second column of a tableaux of size $2k-1$. Conversely, each tableaux of size $2k-1$ yields exactly two tableaux of size $2k$, one where both columns have even length and one where both have odd length.  

By \cite[Proposition 3]{simionschmidt}, for any $n$ we have 
\[ \vert \ivn(123) \vert = \binom{n}{\lfloor \frac{n}{2} \rfloor}\,.\]

Since \[\frac{1}{\sqrt{1-4t^2}} = \sum_{n=0} \binom{2n}{n}t^{2n}\,,\]
we have \[\frac{1}{2t^2\sqrt{1-4t^2}} = \sum_{n=0} \binom{2n+1}{n} t^{2n} + \frac{1}{2t^2}\,,\] and 
\[
\frac{1-\sqrt{1-4t^2}}{2t^2\sqrt{1-4t^2}} = \sum_{n=0} \binom{2n+1}{n}t^{2n}\,.\]
Multiplying this term by the factor $(xt+(1+x^2)t^2)$, and adding 1, gives the desired generating function.
\qed

\section*{Appendix}

In this appendix we restate various results from Flajolet and Sedgewick so they are more readily available for the reader.  

\begin{theorem}[Theorem IX.1 \cite{FlSe09}]\label{IX.1}
Let $\Omega$ be an arbitrary set contained in the unit disc and having at least one accumulation point in the interior of the disc.  Assume that the probability generating function $p_n(u) = \sum_{k\geq 0} p_{n,k}u^k$ and $q(u)=\sum_{k\geq 0} q_ku^k$ are such that there is convergence,
\[\lim_{n\to\infty} p_n(u)=q(u),\]
pointwise for each $u$ in $\Omega$.  Then a discrete limit law holds in the sense that, for each $k$,
\[ \lim_{n\to\infty} p_{n,k}=q_k \quad \textrm{and} \quad \lim_{n\to\infty} \sum_{j\leq k} p_{n,j}= \sum_{j\leq k} q_j.\]
\end{theorem}

The following is essentially \cite[Theorem IX.9]{FlSe09}, rephrased to be more easily applied in the present context and including some of the intermediate results from the proof.  If $f$ is twice differentiable near $1$ we define
\[ \mathfrak{m}(f)= \frac{f'(1)}{f(1)} \quad \textrm{and}\quad \mathfrak{v}(f) = \frac{f''(1)}{f(1)} +\frac{f'(1)}{f(1)} - \left(\frac{f'(1)}{f(1)}\right)^2.\]

\begin{theorem}[Theorem IX.9 \cite{FlSe09}]\label{IX.9}
Let $F(u,z)$ be a function that is bivariate analytic at $(u,z)=(0,0)$ and has non-negative coefficients.  Assume that $F(1,z)$ is meromorphic in $z\leq r$ with only a simple pole at $z=\rho$ for some positive $\rho<r$.  Assume also the following conditions.
\begin{enumerate}
\item Meromorphic perturbation: there exists some $\epsilon>0$ and $r>\rho$ such that in the domain $\mathcal{D} = \{|u-1|<\epsilon\}\times \{|z|\leq r\}$ the function $F(u,z)$ admits the representation
\[ F(u,z) = \frac{B(u,z)}{C(u,z)},\]
where $B(u,z)$ and $C(u,z)$ are analytic for $(u,z)\in \mathcal{D}$, with $B(1,\rho)\neq 0$. (Thus $\rho$ is a simple zero of $C(1,z)$.)
\item Non-degeneracy: one has $\partial_zC(1,\rho)\cdot\partial_u(1,\rho)\neq 0$, ensuring the existence of a non-constant $\rho(u)$ analytic at $u=1$, such that $C(u,\rho(u))=0$ and $\rho(1)=\rho$.
\item Variability: one has 
\[\mathfrak{v}\left(\frac{\rho(1)}{\rho(u)}\right)\neq 0.\]
\end{enumerate}
Let $X_n$ be a random variable with probability generating function
\[ p_n(u) = \frac{[z^n]F(u,z)}{[z^n]F(1,z)},\]
and let $Z$ be a standard normal random variable.  Then, for all $x \in \R$,
\[ \P\left( \frac{X_n-\E(X_n)}{\sqrt{Var(X_n)}} \leq x \right) = \P(Z\leq x)+O\left(\frac{1}{\sqrt{n}}\right).\]
Furthermore, 
\[ \E(X_n) = \mathfrak{m}\left(\frac{\rho(1)}{\rho(u)}\right)n + O(1) \quad \textrm{and}\quad Var(X_n) = \mathfrak{v}\left(\frac{\rho(1)}{\rho(u)}\right)n+O(1).\]
\end{theorem}

For the next result, we need the following definition.

\begin{definition}[Definition VI.I \cite{FlSe09}]
Given two number $\phi$ and $R$ with $R>1$ and $0<\phi<\pi/2$, the open domain $\Delta(\phi,R)$ is defined as
\[\Delta(\phi,R) = \{z\ | \ |z|<R, \ z\neq 1, \ |\arg(z-1)|>\phi\}.\]
For a complex number $\zeta$ a domain $D$ is a $\Delta$-domain at $\zeta$ if there exist $\phi$ and $R$ such that $D=\zeta \Delta(\phi,R)$.   A function is $\Delta$-analytic if it is analytic on a $\Delta$-domain.
\end{definition}  

\begin{corollary}[Corollary VI.1 \cite{FlSe09}]\label{VI.1}
Assume that $f(z)$ is $\Delta$-analytic and
\[ f(z)\sim (1-z)^{-\alpha}, \qquad \textrm{as } z\to 1, \quad z\in \Delta,\]
with $\alpha\notin \{0,-1,-2,\dots\}$.  Then the coefficients of $f$ satisfy
\[ [z^n]f(z) \sim \frac{n^{\alpha-1}}{\Gamma(\alpha)}.\]
\end{corollary}

\bibliographystyle{alpha}

\end{document}